\documentclass[12pt]{amsart}
\usepackage{amssymb,amsmath}
\usepackage{latexsym,color}
\usepackage[mathscr]{eucal}

\newcommand{\defn}[1]{Def.{\ref{#1}}}

\newcommand{\Places}{\mathcal{P}}
\newcommand{\PLACES}[1]{(#1)_{\textit{places}}}
\newcommand{\NODES}[1]{(#1)_{\textit{nodes}}}
\newcommand{\PN}[1]{\mathcal{N}(#1)}
\newcommand{\REDS}{\mathcal{F}}
\newcommand{\POWNODES}{\mathcal{Q}}
\newcommand{\infinite}[1]{\card{#1}\geq\aleph_0}

\newcommand{\COMMENT}[1]{}

\newcommand{\card}[1]{{\left|{#1}\right|}}

\newcommand{\defAs}{\stackrel{_{\tiny {\rm def}}}{_{\normalsize =}}}

\newcommand{\iffAs}
 {\;\mbox{iff}\raisebox{-0.5 ex}[0 ex][0 ex]{\tiny
Def}\:}

\newcommand{\opar}{\mbox{\sf(}\,}
\newcommand{\cpar}{\,\mbox{\sf)}}

\newcommand{\Not}{{\bf\neg}}
\renewcommand{\And}{\wedge}

\newcommand{\biimplies}{{\bf\leftrightarrow}}

\newcommand{\anonymous}{\_} 

\newcommand{\Pow}{\wp}
\newcommand{\pow}[1]{\wp({#1})}

\newcommand{\powAst}{\wp^{\ast}}
\newcommand{\powast}[1]{\powAst({#1})}

\newcommand{\TARGETS}{T}
\newcommand{\Targets}[1]{T({#1})}

\def\qed {{\unskip\nobreak\hfil\penalty50
\hskip2em\hbox{}\nobreak\hfil \rule{2mm}{2mm}
\parfillskip=0pt \finalhyphendemerits=0 \par \medskip}}

\newcommand{\seq}[1]{{\mathcal {A}} _{#1}}

\newenvironment{des}{
\begin{list}
 {$\bullet$}
 {\topsep = 0 mm
\labelwidth = 3 mm
\labelsep = 2 mm
\parsep = 0.4 mm
\itemsep = \parskip
\leftmargin = 5 mm}
}{\end{list}}

\newtheorem{mydef}{Definition} 
\newtheorem{mylemma}[mydef]{Lemma} 
\newtheorem{mytheorem}[mydef]{Theorem} 
\newtheorem{mycorollary}[mydef]{Corollary} 
\newtheorem{myremark}[mydef]{Remark} 
\newtheorem{myproblem}[mydef]{Problem} 

\def \au {\rm}
\def \ti {\it}
\def \jou {\rm}

\def \no#1#2#3 {{\bf #1} (#3), #2.}
\def \eds#1#2#3 {#1, #2, #3.}

\begin{document}
\title[A generalized small model property]
{A generalized small model property\\ 
for languages which force the infinity\\
}

\author[P. Ursino]{Pietro Ursino}

\address{Dipartimento di Matematica
\newline\indent
Universit\`a di Catania
\newline\indent
Viale Andrea Doria 6
\newline\indent
I-95125 Catania, Italy}
\email{ursino@dmi.unict.it}

\begin{abstract}
This paper deals with
formulas of set theory which force the infinity.
For such formulas, we provide a technique to infer 
satisfiability from a finite assignment.
\end{abstract}

\maketitle


\section{Introduction}

In 1970 Jacob T.\ Schwartz launched the computable set theory
longterm project \cite{Sch}, which aimed to merge set theory and
theoretical computer science with reciprocal benefits.
Since then, 
this research field revealed its pure combinatorial behavior.

Ten years later, M.\ Breban (cf.\ \cite{BreFe}) made an attempt to solve 
the decidability problem for the language consisting of
the conjuctions of literals of
the following forms:
\begin{equation}
\begin{array}{llll}
    v=w, & v\neq w, & v=\emptyset, & v=u \cup w,\\
    v=u \cap w, & v=u \setminus w, & v \subseteq u, & v
\not\subseteq u,\\
    v\in w, & v\notin w, & v=\pow{w}, &
    v=\{w_{0},w_{1},\dots,w_{H}\},
\end{array} \label{literals}\tag{$\mathbf{\dagger}$}
\end{equation}
Breban was able to solve the problem allowing at most one occurrence 
of the powerset operator.
Indeed, this unquantified language, known as MLSSP (i.e.,\
\emph{Multi-Level Syllogistic with Singleton and Powerset
operators}),
shows how drastically the complexity of combinatorics increases,
as one enriches the language
with new strong set constructors. 
In \cite{Fe} Ferro solved the problem with two occurrences
of the powerset constructor;
whereas Cantone (see \cite{C1}), 
exploiting a more sophisticated approach, solved the whole 
decidability problem for MLSSP,
without any restriction on the number of occurrences. 
However, any attempt to use the same
simple combinatorial approach to lengthen the list of set constructors
(in a non trivial way), crashed against the
fact that such languages build formulas which force 
any model to be infinite.
Therefore, one of the main goals in solving advanced decidability problems
is to find a way to overcome the impossibility 
to find finite models not exceeding a fixed size.

Recently (see \cite{CU97}),
the use of formative processes as a \emph{history} of a set assignment
gave a new perspective to solve this kind of problems.
Indeed, it makes use of the history (or trace) of the model
to obtain new information in order to decrease 
the size of the model up to a suitable one. 
This observation
motivated our interest to the study of a small model property for 
languages which contain MLSSP.

In \cite{Cantone-Omodeo-Ursino02}, we 
discovered the small model property
for MLSSP, and,
by means of this result, we built a satisfiability decision algorithm.

If we add to MLSSP particular set constructors, the small model property fails to hold.
A rather explicit example is the finiteness operator $\mathit{Finite}(x)$
(meaning that the cardinality of the set designated by $x$ is
smaller than $\aleph_0$). 
Of course, since we admit negation
among propositional connectives, we must also take into account
literals of the form $\Not\,\mathit{Finite}(x)$. 
Thus MLSSP, extended with
the monadic relator $\mathit{Finite}$, ``forces the infinity"
(informally, a language forces the infinity whenever
has inside formulas whose models must be of infinite size).
The same happens allowing 
the unitary union operator $\bigcup (x)$. 
As consequence, languages which allow the use of 
this type
of operators cannot satisfy the small model property.
This gave us the suggestion to focus on the structure of infinite models
(in particular, to their combinatorial features).
Hence we formulate

\begin{myproblem}
Which combinatorial properties two assignments 
have to share, in order to satisfy the same 
{\rm MLSSP}-like literals?
\end{myproblem}

Corollary \ref{MsatisfiesII} below gives a satisfactory 
answer to this question. In Corollary~\ref{upwards},
we provide an analogous result, but referred to the 
formative processes of the assignments.

These two corollaries are the tool to prove how
a finite assignment can be equipped with a special structure 
that allows to increase some variables,
without affecting the validity of the formula.
We agree to denote such  variables as \emph{potential infinite variables},
and we find a condition for this property to hold.
The above results allow us to investigate 

\begin{myproblem}
Even if a language forces the infinity, is it still possible, 
for any satisfiable formula, to 
exhibit a finite assignment that
witnesses this satisfiability or, in other terms, to show a finite representation of an infinite model?
\end{myproblem}

This kind of property of languages is here introduced as 
\emph{witness-small model property}. 

Theorem~\ref{pumping1} 
demonstrates how a combinatorial property of a finite assignment to a formula of MLSSPF can
witness the satisfiability of literals which require an infinite assignment.
More generally, this paper shows how
the formative processes can be used in order to prove the witness-small model property 
in some cases.
This result 
leads to the solution of 
some open problems, such as the decidability
of languages which allow the use of the above-cited set constructors, namely,
MLSSP extended with
the monadic relator $\mathit{Finite}$ (the so-called MLSSPF) \cite{COU2}, and 
MLSSP extended with
the monadic operator $\bigcup (x)$ (known as MLSSPU) \cite{CU}.

Our method is based on a specific analysis, both of the model and of a
formative process which generates it. 
An detailed treatment of the general features of computable
set theory can be found in \cite{CFO} and \cite{COP}.

\section{Basic notations and background}

For the reader's convenience, we provide in this section
brief description of the standard
tools used in set computable theory.
For usual set theoretic notion we refer to any textbook of the field
(see \cite{JEC}, for example), instead
a complete survey of the specific notions mentioned in the sequel
may be found in \cite[\S 2]{Cantone-Omodeo-Ursino02}.

\subsection{Assignments and models}
Fix allowed forms for literals. A propositional combinations of 
literals of such forms is said a
\emph{formula}.
It is customary to denote \emph{language of set theory}
the family of 
all formulas built with assigned forms of literals.
Assume $\Phi$ is a formula and
let ${\mathcal M}\in\{\:\mbox{\rm sets}\:\}^{\mathcal
X}$ be a set-valued assignment 
defined on the collection ${\mathcal X}_{\Phi}$
of variables in $\Phi$. If ${\mathcal M}$ satisfies all the literals, 
it is said to be
a \emph{model} for $\Phi$.
A model is \emph{rank-bounded} by $k$ if the rank
of any set involved 
in the assignment does not exceed $k$.

\begin{mydef}\label{}\rm 
A language satisfies the \emph{small model property} if there exists
a computable natural function $f$ such that for any given
formula $\Phi$ of that language
and any model 
${\mathcal M}$ of $\Phi$ there is a finite model ${\mathcal M}'$
rank-bounded by $f(\vert {\mathcal X}_{\Phi}\vert)$.
\end{mydef} 

Assume that $\Phi$ is a formula of a language, and 
${\mathcal A}$ is a set assignment to its variables ${\mathcal X}_{\Phi}$.
We say that ${\mathcal A}$ \emph{witnesses the satisfiability of}
$\Phi$ (even if ${\mathcal A}$ is not a model for $\Phi$),
provided that the structure of ${\mathcal A}$ allows to infer the satisfiability of 
$\Phi$.
A formula of set theory \emph{forces the infinity} if it possesses
a variable $x$
such that, for any model ${\mathcal M}$ which satisfies the formula,
${\mathcal M}(x)$ is of infinite size.
From this point of view, a formula which forces the infinity cannot
have a finite model, but it could have a finite assignment which witnesses its
satisfiability. Hence the following
definition makes sense:

\begin{mydef}\label{}\rm 
A language satisfies the \emph{witness-small model property} 
if there exists
a computable natural function $f$ such that for any given
formula $\Phi$ of that language
and any model 
${\mathcal M}$ of $\Phi$
there exists a finite assignment ${\mathcal A}$
rank-bounded by $f(\vert {\mathcal X}_{\Phi}\vert)$ which witnesses the
satisfiability of $\Phi$. 
\end{mydef} 

\subsection{Transitive partitions and syllogistic boards}

\begin{mydef}\label{defPart}\rm
A family $\Sigma$ of pairwise nonempty disjoint sets
is called a \emph{partition} (of
$\bigcup\Sigma$). Its members are the \emph{blocks} of
$\Sigma$.
The set
$\varsigma_{\Sigma}\defAs\pow{\bigcup\Sigma}\setminus\bigcup\Sigma$
(often denoted simply by $\varsigma$) will occasionally be
treated as a block of the partition too. In this case, it is called the
\emph{outer block} of $\Sigma$.
\end{mydef}

As is well known, the function
$$\Sigma\stackrel{\sim_{\Sigma}}{\longmapsto}\{[X,Y]\,
\mbox{\tt|}\:(\exists\,b\in\Sigma)( X\in b\And Y\in b)\}
$$ 
establishes
a one-to-one correspondence between the partitions of a given set $S$
and the equivalence relations on $S$.

A useful relation $\sqsubseteq$ on $\pow{\pow{ S}}$ is defined by
setting
$${\mathcal B}\sqsubseteq{\mathcal
A}\;\iffAs\;(\forall\,a\in{\mathcal A})(\exists\,B\subseteq{\mathcal
B})\:a=\bigcup B\:.$$
The relation ${\mathcal B}\sqsubseteq{\mathcal A}$ reads
`` $\mathcal B$ is {\em finer} than $\mathcal A$ ", or
`` $\mathcal A$ is {\em coarser} than $\mathcal B$". This obviously is a
{\em preorder} relation that, when restricted to the set
$\varpi( S )$ of all partitions of
$S$, $\sqsubseteq$, becomes a {\em partial ordering}.

\begin{mydef}\label{transitivePartition}\rm
A partition $\Sigma$ is
said to be \emph{transitive} if $\bigcup\Sigma$ is transitive.
\end{mydef}


We consider a finite set
$\Places$, whose elements are called \emph{places} and whose
subsets are called \emph{nodes}. Places and nodes will be the
vertices of a directed bipartite graph $\mathcal{G}$ of a special
kind, called a \emph{$\Places$-board}. The edges issuing from
each place $q$ are, mandatorily, all pairs $q,B$ such that $q\in
B\subseteq \Places$. The remaining edges of $\mathcal{G}$ must
lead from nodes to places. Hence, $\mathcal{G}$ is fully
characterized by the so called \emph{target function}
$$\TARGETS\:\in\:\pow{\Places}^{\pow{\Places}},$$
associating with each node $A$ the set of all places $t$ such
that $\langle A,t\rangle$ is an edge of $\mathcal{G}$. The
elements of $\Targets{A}$ are called the \emph{targets} of
$A$. We will usually represent $\mathcal{G}$ simply by
$\TARGETS$.

Places and nodes of a $\Places$-board are meant to represent the
blocks $\sigma$, and the subsets $\Gamma$ (or, quite often, their
unionsets $\bigcup\Gamma$), of a transitive partition $\Sigma$,
respectively. Moreover, in this case, 
there is a quite natural way to define the above-mentioned 
directed bipartite graph structure.

For our convenience we define the further operator
$$A\ni\in B=_{def}A\cap B\ne\emptyset.$$
For any set $X$, we put
$$\powast{ X }\defAs\{\,Y\,\mbox{\tt|}\:Y\subseteq \bigcup X\And
(\forall\,z\in X)(\:z\ni\in Y\:)\,\}\:,$$
that is, the elements of the family $\powast{X}$ are all the
sets $Y$ that can be obtained by extracting from each $z\in X$ a nonnull
$W_{z}\subseteq z$, so forming $Y=\bigcup_{z\in X}W_{z}$.

\begin{mydef}\rm
A transitive partition $\Sigma$ is said to
\emph{comply with} $\mathcal{G}$ via $q\mapsto
q^{(\bullet)}$, where $\mathcal{G}$ is
$\Places$-board,
$q\mapsto q^{(\bullet)}$ belongs to
$\Sigma^{\Places}$
and $\TARGETS (A)=\{q\mid\powast{A^{(\bullet)}}
\ni\in q^{(\bullet)}\}$, if the function $\TARGETS$ satisfies all the
properties required by $\mathcal{G}$, as indicated above (in
particular, this requires $q\mapsto q^{(\bullet)}$ to be
injective).
\end{mydef}

Any such board is said to be \emph{induced} by $\Sigma$
(for short, a \emph{$\Sigma$-board}). 
We denote a transitive $\Sigma$-board by a couple $(\Sigma,\mathcal{G})$, 
where $\Sigma$ is a transitive partition and
$\mathcal{G}$ is the induced $\Places$-board.

For the purposes of this paper, some additional structure must be
superimposed on $\Places$-boards:

\begin{mydef}\rm
A $\Places$-board $\mathcal{G}=(\TARGETS,\REDS,\POWNODES)$ 
is said to be \emph{colored} when it has
\begin{itemize}
\item a designated set $\REDS$ of places, 
\item a designated set $\POWNODES$ of nodes, such that
$D\in\POWNODES$ holds whenever $D\subseteq B\in\POWNODES$ \\ (in
short, $\bigcup\Pow[\POWNODES]\subseteq\POWNODES$ ), and
\item a target function $T$.
\end{itemize}
The places in $\REDS$ are said to be \emph{red}, the ones in
$\Places\setminus\REDS$ are said to be \emph{green}; the nodes
in $\POWNODES$ are called \emph{$\Pow$-nodes}.
A node is red if all places in it are red, and
green otherwise; a list of vertices is 
green if all vertices lying on it are green.
\end{mydef}

\begin{mydef}\label{SimBoard}\rm 
Let ${\mathcal G}$ be a colored transitive \emph{$\Sigma$-board}. Then
$\widehat{\Sigma}$ is said to \emph{simulates}
$(\Sigma,{\mathcal G})$ \emph{upwards}, 
when there is a bijection
$\beta\in\widehat{\Sigma}^{\Sigma}$ such that
\begin{itemize}
\item $\widehat{\Sigma}$ $\in$-simulates $\Sigma$ via
$\beta$. That is, $\bigcup\beta[X]\in\bigcup\beta[Y]$ if and only if
$\bigcup X\in\bigcup Y$, for $X,Y\subseteq\Sigma$;
\item $\widehat{\Sigma}$ $\Pow$-simulates $\Sigma$ via
$\beta$. That is, $\bigcup\beta[X]=\pow{\bigcup\beta[Y]}$ if 
$\bigcup
X=\pow{\bigcup Y}$, for $Y\in\POWNODES$ $X,Y\subseteq\Sigma$.
\item $\widehat{\Sigma}$ $Red$-simulates $\Sigma$ via
$\beta$. That is,
if $\sigma\in\REDS$, then $\card{\beta(\sigma)}=\card{\sigma}$;

\end{itemize}
\end{mydef}

As far as the Boolean constructs
$\emptyset,\cap,\setminus,\cup,=,\neq,\subseteq, \not\subseteq$
are concerned, all relevant information about a family of sets is
conveyed by the following structure:

\begin{mydef}\rm
Given a family $\mathcal F$, the \emph{Venn partition} of $\mathcal F$
is the coarsest partition $\Sigma$ of ~$\bigcup{\mathcal F}$ which
fulfill the condition
$$(\forall\,x\in{\mathcal
F})(\forall\,p\in\Sigma)(\:p\ni\in x\;\rightarrow\;p\subseteq
x\:).$$
\end{mydef}

Assume that $\Phi$ is a collection of literals which have one of the forms
$($\ref{literals}$)$, and let
${\mathcal M}\in\{\:\mbox{\rm sets}\:\}^{{\mathcal X}_{\Phi}}$
be a set-valued assignment defined on the collection 
${\mathcal X}_{\Phi}$
of variables in $\Phi$.
We denote by $\Sigma_{{\mathcal X}_{\Phi}}$
the Venn partition of the set ${\mathcal M}[{\mathcal X}_{\Phi}]$, and by $\Im_{{\mathcal M}}$
the function $\Im_{{\mathcal M}}\in\pow{\Sigma_{{\mathcal X}_{\Phi}}}^{{\mathcal X}_{\Phi}}$ such that
${\mathcal M}(v)= \bigcup\Im_{{\mathcal M}}(v)$ holds for every $v$ in ${\mathcal X}_{\Phi}$.

\begin{myremark}\rm
Observe that any formula $\Phi$ with variables ${\mathcal X}_{\Phi}$
of a language resulting from an extension of Multi Level Syllogistic can be modified,
without affecting its satisfiability, in such a way any model ${\mathcal M}$ generates a transitive
$\Sigma_{{\mathcal X}_{\Phi}}$ \cite[pp.195-196]{C1}.
Because of that, from now on we shall assume that
$\Sigma_{{\mathcal X}_{\Phi}}$ is transitive, for any model ${\mathcal M}$ 
of a formula $\Phi$
with variables ${\mathcal X}_{\Phi}$.
\end {myremark}

Whenever literals as $v=\pow{w}$ and $Finite(v)$ 
appear in $\Phi$,  $\Sigma_{{\mathcal X}_{\Phi}}$ can
be naturally transformed into a colored \emph{$\Sigma_{{\mathcal X}_{\Phi}}$-board} 
${\mathcal G}=(\TARGETS,\REDS,\POWNODES)$ 
(i.e., the $\Sigma$-board ${\mathcal G}$ 
induced by $\Sigma_{{\mathcal X}_{\Phi}}$), in the following way.
\begin{itemize}
\item[(a)] $\REDS=\bigcup\{\Im(v)\mid \mbox{ for all literals of the 
form $v=\{w_{1},\dots,w_{H}\}$  and $Finite(v)$ in $\Phi$}\}$; 
\item[(b)] $\POWNODES$ is equal to the minimal collection
of nodes such that
\begin{itemize}
\item $\Im(u)\in\POWNODES$ for all literals of the form $u=\pow{w}$ in $\Phi$, and 
\item $\bigcup\Pow[\POWNODES]\subseteq\POWNODES$.
\end{itemize}
\end{itemize}
In the above case we refer to such a $\Sigma_{{\mathcal X}_{\Phi}}$-board as 

\emph{ the canonical board of the assignment ${\mathcal M}$ to the $MLSSPF$ formula $\Phi$}.
\begin{mylemma}
\label{MsatisfiesI}
Consider a formula $\Phi\in {\rm MLSSPF}$,
a set-valued assignment ${\mathcal M}\in\{\:\mbox{\rm sets}\:\}^{{\mathcal X}_{\Phi}}$ 
defined on the collection ${\mathcal X}_{\Phi}$
of variables in $\Phi$, together with the 
colored transitive $\Sigma_{{\mathcal X}_{\Phi}}$-board 
${\mathcal G}=(\TARGETS,\REDS,\POWNODES)$.
Define $\Phi^-$ as the formula $\Phi$ without literals of the type $Finite(x)$
or $\neg Finite(x)$.
Moreover, let be $\widehat{\Sigma}$ a partition and $\beta$  a
bijection between $\Sigma_{{\mathcal X}_{\Phi}}$ and $\widehat{\Sigma}$ such that 
$\widehat{\Sigma}$ simulates 
$(\Sigma,{\mathcal G})$ upwards via
$\beta$, and let ${\mathcal M}'(v)=\bigcup\beta[\Im_{{\mathcal M}}(v)]$.
Then, for every literal in $\Phi^-$, the following conditions are
fulfilled:
\begin{des}
\item if the literal is satisfied by $\mathcal M$, then it is satisfied
by ${\mathcal M}'$ too;
\item if the literal is satisfied by ${\mathcal M}'$, and does not
involve $\Pow$~ or the construct $\{\anonymous,\dots,\anonymous\}$,
then it is satisfied by ${\mathcal M}$ too;
\item if the literal $Finite(x)$ appears in $\Phi$ and is satisfied by 
$\mathcal M$, then it is satisfied
by ${\mathcal M}'$ too.
\end{des}
\end{mylemma}

\begin{proof}
The thesis can be recast as follows. For $u,v,w$ and $w_{i}$ in
${\mathcal X}_{\Phi}$, the following conditions hold for all literals in $\Phi$:
\begin{itemize}
\item[(1)] $\bigcup\Im(v)\;\Re\;\bigcup\Im(w)$ iff
$\bigcup\beta[\Im(v)]\;\Re\;\bigcup\beta[\Im(w)]$,
 for $\Re$ in $\{\:=\,,\:\in\,,\:\subseteq\:\}$;
\item[(2)] $\bigcup\Im(v)=\bigcup\Im(u)\:\star\:\bigcup\Im(w)$ iff
 $\bigcup\beta[\Im(v)]=\bigcup\beta[\Im(u)]\:\star\:
 \bigcup\beta[\Im(w)]$,
 for $\star$ in $\{\:\cap\,,\:\setminus\,,\:\cup\:\}$, and 
$\bigcup\Im(v)=\emptyset$ iff
$\bigcup\beta[\Im(v)]=\emptyset$;
\item[(3)] if $\bigcup\Im(v)=\pow{\bigcup\Im(w)}$, then
$\bigcup\beta[\Im(v)]=\pow{\bigcup\beta[\Im(w)]}$;
\item[(4)] if 
$\bigcup\Im(v)=\{\,\bigcup\Im(w_{1}),\dots,\bigcup\Im(w_{H})\,\}$,
then\\ $\bigcup\beta[\Im(v)]=\{\,\bigcup\beta[\Im(w_{1})],\dots,
\bigcup\beta[\Im(w_{H})]\,\}$.
\item[(5)] if $Finite(v)$ appears in $\Phi$ then 
$\vert\bigcup\Im(v)\vert=\vert\bigcup\beta[\Im(v)]\vert$
\end{itemize}
Property $(1)_{\in}$ (here $\Re$ is meant to be $\in$) follows from $\in$-simulates in
\defn{SimBoard}. $(3)$ follows from the assumption $\Im(v)\in\POWNODES$ and the
notion of $\Pow$-simulates given in the same definition. 
Condition $(5)$ plainly follows from definition of $Red$-simulates.

We are left to prove that $(4)$ hold. 
Observe that $\Im(v)\subseteq\REDS$, then consider $\Im(v)$
as the set $X$ and $Y_i$ as the sets $\Im(w_{i})$. Hence we can assume that
$\bigcup X=\{\bigcup Y_1,\dots, \bigcup Y_L\}$, $X\subseteq\REDS$, and
$Y_1,\dots, Y_L$ are distinct. We must check that
$\bigcup\beta[X]=\{\bigcup\beta[Y_1],\dots,$ $\bigcup\beta[Y_L]\}$.
Since $\widehat{\Sigma}$ $Red$-simulates
$(\Sigma,{\mathcal G})$ and $X\subseteq\REDS$,
and $\card{\beta(\sigma)}=\card{\sigma}$ for
each $\sigma\in X$, the desired conclusion easily follows.
Indeed, by property (1) of Def.~\ref{SimBoard}, 
$\bigcup\beta[Y_i]\in\beta(\sigma)$
if and only if
$\bigcup Y_i\in\sigma$, and $\beta[Y_1],\dots,\beta[Y_L]$
(and, accordingly, $\bigcup\beta[Y_1],\dots,\bigcup\beta[Y_L]$)
are pairwise distinct.

The proofs of remaining bi-implications go exactly as in \cite[Lemma 10.1]
{Cantone-Omodeo-Ursino02}
\end{proof}

\begin{mydef} 
\label{defImitation}\rm
Consider a colored \emph{$\Sigma$-board} ${\mathcal G}=(\TARGETS,\REDS,\POWNODES)$
A partition $\widehat{\Sigma}$ is said to \emph{imitate}
$(\Sigma,{\mathcal G})$ when there is a bijection $\beta\in
\widehat{\Sigma}^{\Sigma}$ such that, for
$\Gamma\subseteq\Sigma$, $\sigma\in\Sigma$,
\begin{itemize}
\item[(1)] $\beta(\sigma)\ni\in\powast{\beta[\Gamma]}$ holds [if
and]
only if $\sigma\ni\in\powast{\Gamma}$;
\item[(2)] $\bigcup\beta[\Gamma]\in\beta(\sigma)$ holds if and only if $\bigcup
\Gamma\in\sigma$;
\item[(3)] if $\Gamma\in\POWNODES$ holds, then
$\powast{\beta[\Gamma]}\subseteq\bigcup\widehat{\Sigma}$;
\item[(4)]  if
$\sigma\in\REDS$ holds, then $\card{\beta(\sigma)}<\aleph_0$.
\end{itemize}
We will say that $\widehat{\Sigma}$
\emph{imitates} $(\Sigma,{\mathcal G})$
\emph{upwards} when the following
additional condition holds, for all $\sigma\in\Sigma$:
\begin{itemize}
\item[(4$'$)] if
$\sigma\in\REDS$, then $\card{\beta(\sigma)}=\card{\sigma}$;
\end{itemize}
\end{mydef}

\begin{mylemma}\label{imitaSimula}
Consider a colored \emph{$\Sigma$-board} ${\mathcal G}=(\TARGETS,\REDS,\POWNODES)$
assume that a transitive partition $\widehat{\Sigma}$ {imitates} $(\Sigma,{\mathcal G})$
 {upwards}
then it {simulates} 
$(\Sigma,{\mathcal G})$ {upwards}.
\end{mylemma}

\begin{proof}
Let $\Sigma$ and $\widehat{\Sigma}$ be transitive partitions, and  
let ${\mathcal G}$ be a colored \emph{$\Places$-board} induced by 
$\Sigma$. Assume that $\widehat{\Sigma}$ {imitates}
$(\Sigma,{\mathcal G})$ {upwards}
via the bijection $\beta\in
(\widehat{\Sigma})^{\Sigma}$.
Finally, let $X, Y \subseteq \Sigma$. 

Then we have: $\bigcup\beta[X]\in\bigcup\beta[Y]$ iff
$(\exists\,\widehat{\sigma}\in\beta[Y])(\bigcup\beta[X]\in\widehat{\sigma})$
iff $(\exists\,\sigma\in
Y)(\bigcup\beta[X]\in\beta(\sigma))$
iff $(\exists\,\sigma\in
Y)(\bigcup X\in\sigma)$ iff $\bigcup X\in\bigcup Y$.

Assuming now that $\bigcup X=\pow{\bigcup Y}$, $Y\in\POWNODES$, let us 
prove that
$\pow{\bigcup\beta[Y]}\subseteq\bigcup \beta[X]$.  Indeed, suppose
$t\subseteq\bigcup\beta[Y]$ and let $\widehat{\Sigma}_{t}$ be the
subset of $\widehat{\Sigma}$ for which
$t\in\powast{\widehat{\Sigma}_{t}}$ (so that
$\widehat{\Sigma}_{t}\subseteq \beta[Y]$, which implies
$\widehat{\Sigma}_{t}\in\POWNODES$ by the hereditarily closedness by inclusion
of $\POWNODES$).  As
$\beta^{-1}[\widehat{\Sigma}_{t}] \subseteq Y$,
it follows that
$\powast{\beta^{-1}[\widehat{\Sigma}_{t}]} \subseteq
\pow{\bigcup Y} =
\bigcup X \subseteq \bigcup \Sigma$.
Therefore, by the fact that 
$\widehat{\Sigma}$ {imitates}
$(\Sigma,{\mathcal G})$ {upwards} and $\widehat{\Sigma}_{t}\in\POWNODES$,
it follows that 
$\powast{\widehat{\Sigma}_{t}} \subseteq \bigcup \widehat{\Sigma}$,
so
that $t \in \bigcup \widehat{\Sigma}$.  Let $\widehat{\sigma}_{t}$
be the block in $\widehat{\Sigma}$ to which $t$ belongs, and let
$\sigma_{t}$ be the block in $\Sigma$ for which
$\beta(\sigma_{t})=\widehat{\sigma}_{t}$.  Then, since
$\powast{\widehat{\Sigma}_{t}}\ni\in\widehat{\sigma}_{t}$, we 
have that
$\powast{\beta^{-1}[\widehat{\Sigma}_{t}]}\ni\in\sigma_{t}$, which
yields $\bigcup X=\pow{\bigcup Y}\supseteq\powast{\beta^{-1}[
\widehat{\Sigma}_{t}]}\ni\in\sigma_{t}$, so that $\bigcup
X\ni\in\sigma_{t}$, $\sigma_{t}\in X$, and hence
$t\in\widehat{\sigma}_{t} \in\beta[X]$, which in turn yields $t \in
\bigcup \beta[X]$.

Next, assuming again 
$\bigcup X=\pow{\bigcup Y}$, let us prove that $\bigcup
\beta[X]\subseteq\pow{\bigcup\beta[Y]}$.  Indeed, for each $t\in
\bigcup\beta[X]$ there is a unique $\sigma_{t}\in X$ such that
$t\in\beta(\sigma_{t})$; moreover, by the transitivity of
$\bigcup\widehat{\Sigma}$, there is a unique $\Gamma\subseteq\Sigma$ for
which $t\in\powast{\beta[\Gamma]}$. Moreover, since
$\powast{\beta[\Gamma]}\ni\in\beta(\sigma_{t})$, we also have that
$\powast{\Gamma}\ni\in\sigma_{t}$.
Thus we can take 
$t'\in\sigma_{t}\cap\powast{\Gamma}$ that, as
$\sigma_{t}\subseteq\bigcup X=\pow{\bigcup Y}$, fulfills
$t'\in\powast{Z}$ for a suitable $Z\subseteq Y$.  In conclusion,
$\Gamma=Z$, and therefore
$t\subseteq\bigcup\beta[\Gamma]=\bigcup\beta[Z]\subseteq\bigcup\beta[Y]$.

\end{proof}

As an immediate consequence, we have

\begin{mycorollary}\label{MsatisfiesII}
Consider a formula $\Phi\in {\rm MLSSPF}$,
a set-valued assignment ${\mathcal M}\in\{\:\mbox{\rm sets}\:\}^{{\mathcal X}_{\Phi}}$ 
defined on the collection ${\mathcal X}_{\Phi}$
of variables in $\Phi$, together with the 
colored transitive {$\Sigma_{{\mathcal X}_{\Phi}}$-board} 
${\mathcal G}=(\TARGETS,\REDS,\POWNODES)$.
Moreover, let $\widehat{\Sigma}$ and $\beta$ be a partition and a
bijection, respectively, 
such that $\widehat{\Sigma}$ {imitates} $(\Sigma,{\mathcal G})$
{upwards} via
$\beta$, and let ${\mathcal M}'(v)=\bigcup\beta[\Im_{{\mathcal M}}(v)]$, 
where $\Im$ is
the function $\Im\in\pow{\Sigma}^{\mathcal X_{\Phi}}$ such that
${\mathcal M}(v)= \bigcup\Im(v)$ holds for every $v$ in $\mathcal X$.
Then, for every literal in $\Phi^-$ and literals of the type
$Finite(x)$, the following conditions are
fulfilled:
\begin{des}
\item if the literal is satisfied by $\mathcal M$, then it is satisfied
by ${\mathcal M}'$ too;
\item if the literal is satisfied by ${\mathcal M}'$, and does not
involve $\Pow$~ or the construct $\{\anonymous,\dots,\anonymous\}$,
then it is satisfied by ${\mathcal M}$ too.
\end{des}
\end{mycorollary}

\subsection{Formative processes}
We now formalize the concept of ``history'' of a model 
by a transfinite construction.
Using the transitivity of any transitive partition,
it is possible to single out a process that
builds it, having the empty partition as starting point.

The following notions are introduced to specify this concept.

\begin{mydef}\label{prolongation}\rm
Let $\Sigma$ and $\Sigma'$ be two partitions, and let $\Gamma
\subseteq \Sigma$. We say that $\Sigma'$
{\em prolongates} $\Sigma$ via $\Gamma$ when the following
conditions hold:
\begin{enumerate}
\item
for all $\sigma\in\Sigma$, there is one and only one
$\sigma'\in\Sigma'$ such that $\sigma\subseteq\sigma'$;
\item
$\bigcup\Sigma'\setminus\bigcup\Sigma\subseteq\powast{\Gamma}$;
\item $\Sigma\neq\Sigma'$.
\end{enumerate}
When just condition (1) is met, possibly without (2) or (3),
we say that $\Sigma'$ {\em extends} $\Sigma$. If
both (1) and (3) hold true, then $\Sigma'$ is said to extend $\Sigma$
{\em properly}.
\end{mydef}

\begin{mydef}{\bf [Coherence requirement]}
\label{cohExtends}\rm
Let $\Gamma$, $\Sigma'$ and $\Sigma''$ be partitions, with
$\Sigma'$ extending $\Gamma$ (typically,
$\Gamma\subseteq\Sigma'$) and $\Sigma''$ extending $\Sigma'$.
Then $\Sigma''$ is said to extend $\Sigma'$ {\em coherently} with
$\Gamma$ if no element of $\bigcup\Sigma''$ belongs to
$\powast{\Gamma}\setminus\bigcup\Sigma'$. 
\end{mydef}

\begin{mydef}\label{formativeProcess}\rm
Let $\xi$ be an ordinal and let
$\left(\{q^{(\mu)}\}_{q\in\Places}\right)_{\mu\leqslant\xi}$ be a
$(\xi+1)$-sequence of functions, all defined on the
same domain $\Places$. Put
$B^{(\mu)}\defAs\{\,q^{(\mu)}\,\mbox{\tt|}\:q\in B\,\}$ for all
$B\subseteq P$, and let
$\Sigma_{\mu}=\Places^{(\mu)}\setminus\{\emptyset\}$, for all
$\mu\leqslant\xi$.

Assume that the following conditions are fulfilled:
\begin{itemize}
\item $q^{(\mu)}\cap p^{(\mu)}=\emptyset$ when $p,q\in P$, $p\neq q$,
and $\mu\leqslant\xi$;
\item $q^{(\nu)}\subseteq q^{(\nu+1)}$ for all $q\in
\Places$ when $\nu<\xi$;
\item
$q^{(\lambda)}=\bigcup_{\nu<\lambda}q^{(\nu)}$ for every $q\in
\Places$ and every limit ordinal $\lambda\leqslant \xi$;
\item
$q^{(0)}=\emptyset$ and $\emptyset\neq q^{(\xi)}$, for all
$q\in\Places$.
\end{itemize}
In particular, $\Sigma_{0}=\emptyset$ and, for every $\mu\leqslant\xi$,
$\Sigma_{\mu}$ is a partition of the subset $\bigcup P^{(\mu)}$
of $\bigcup P^{(\xi)}$.

Assume moreover that to each $\nu<\xi$ corresponds
$\Gamma_{\nu}\subseteq\Sigma_{\nu}$ such that
\begin{itemize}
\item $\Sigma_{\nu+1}$ prolongates $\Sigma_{\nu}$ via
$\Gamma_{\nu}$ (cf.\ \defn{prolongation});
\item $\Sigma_{\xi}$ extends $\Sigma_{\nu+1}$ coherently with
$\Gamma_{\nu}$ (cf.\ \defn{cohExtends}).
\end{itemize}

Then the sequence
$\left(\{q^{(\mu)}\}_{q\in\Places}\right)_{\mu\leqslant\xi}$
(occasionally, $\big(\Sigma_\mu\big)_{\mu\leqslant\xi}$) is
called a ({\em strong}) {\em formative process} for
$\Sigma_{\xi}$.
Furthermore, the $\xi$-sequences $(A_{\nu})_{\nu<\xi}$ and
$(A_{\nu},T_{\nu})_{\nu<\xi}$, with $ A_{\nu}, T_{\nu}\subseteq
\Places$, satisfying for each $\nu$ the conditions
\begin{itemize}
\item $A_{\nu}^{(\nu)}=\Gamma_{\nu}$,
\item $\{\,q^{(\nu+1)}\setminus q^{(\nu)}\,\mbox{\tt|}\:q\in
T_{\nu}\}$ is a partition of $\bigcup\Sigma_{\nu+1}\setminus
\bigcup\Sigma_{\nu}$ ~
($=\opar\powast{\Gamma_{\nu}}\setminus\bigcup\Sigma_{\nu}\cpar
\cap\bigcup\Sigma_{\nu+1}$)
\end{itemize}
are called the {\em trace} of the
formative process, and a {\em history} of $\Sigma_{\xi}$,
respectively.

A {\em weak formative process} is like a formative process, except
that the coherence requirement is withdrawn from the definition.
A {\em weak trace} is defined similarly.
\end{mydef}

In the sequel it will be helpful the following simplified notation.

\begin{mydef}\rm
Let $(\{q^{(\mu)}\}_{q\in\Places})_{\mu\leqslant\xi}$ be a weak
formative process.  Then, for $q\in \Places$, $B\subseteq
\Places$ and $\nu<\xi$, we set
$$
q^{(\bullet)}\defAs q^{(\xi)},\qquad
B^{(\bullet)}\defAs B^{(\xi)},\qquad
\Delta^{(\nu)}(q) \defAs q^{(\nu+1)}\setminus\bigcup\Places^{(\nu)}.
$$
\end{mydef}

If we take, along with a colored $\Places$-board
$(\TARGETS,\REDS,\POWNODES)$, a bijection $q\mapsto q^{(\bullet)}$
from the places $\Places$ to the final partition $\Sigma_\xi$ of
a formative process, and if moreover $\Sigma_\xi$ complies with
$\TARGETS,\REDS,\POWNODES$, we get what we call a
{\em colored $\Places$-process}: namely, the quintuple
$((\{q^{(\mu)}\}_{q\in\Places})_{\mu\leqslant\xi},(\bullet),
\TARGETS,\REDS,\POWNODES)$.

\begin{mydef}\rm
$e\in\bigcup\Places^{(\bullet)}$ is said to be \em{unused}
at $\mu\leqslant\xi$ if $e\notin
\bigcup\bigcup\Places^{(\mu)}$, i.e., if $e\notin z$ for any
$q\in\Places$ and any $z\in q^{(\mu)}$.
\end{mydef}

\begin{mydef}\rm
An $e\in\bigcup\Places^{(\bullet)}$ is said to be \textsc{new}
at $\mu\leqslant\xi$ if $e\in\Delta^{(\mu)}(q)$ for some $q\in\Places$.
\end{mydef}

Obviously a new element is, in particular, unused.

\begin{mylemma}\label{unused}
If $b$ is a set made of unused elements only, the same is
$\powast{\{b\}\cup A}$.
\end{mylemma}

\subsection{Grand events and local trash}

We begin with the following easy remark.
The block at place $s$ belonging to a
$\Pow$-node $A$ cannot become infinite during a colored process,
unless $A$ has a green place among its targets. To see that,
assume that $s\in A\in\POWNODES$ and $\infinite{s^{(\bullet)}}$.
Consequently, $\card{\powast{\bigcup A^{(\bullet)}}}>\aleph_0$ and
$\powast{\bigcup
A^{(\bullet)}}\subseteq\bigcup\Places^{(\bullet)}$. Hence
there must be a place $g$ such that $\card{\powast{\bigcup
A^{(\bullet)}}\cap g^{(\bullet)}}>\aleph_0$, since
$\card{\Places^{(\bullet)}}=\card{\Places}<\aleph_0$. This
obviously implies that $g\in\Targets{A}\setminus\REDS$.

In light of generalizing the above remark, recalling the notion
of grand move, and noticing that such an event occurs, in a colored
process, at most once for each node $A$, we give the following
definition of \emph{grand event} $GE(A)$ associated with $A$.

\begin{mydef}\rm
For every node $A$ and every $\nu$ such that
$0\leqslant{\nu}<\xi$
$$
GE(A)\;\defAs\;\left\{\begin{array}{ll}
\mbox{\rm the ordinal $\nu$ for which $\bigcup
A^{(\bullet)}\in\bigcup\Places^{(\nu+1)}\setminus\bigcup\Places^{(\nu)}$,}
& \mbox{\rm if any exists}, \\
\mbox{\rm the length $\xi$ of the process}, & \mbox{\rm otherwise}.
\end{array}\right.
$$
Moreover, for any given collection ${\mathcal A}$ of nodes, we put
$$
GE({\mathcal A})\;\defAs\; \min \{GE(A)\mid A \in {\mathcal A}\}\,.
$$
\end{mydef}

Notice that this Definition implies that for any node $A$ and any
$\nu$ such that $0\leqslant\nu<\xi$,
$$\begin{array}{lcl} \nu\leqslant GE(A) & \biimplies &
\bigcup A^{(\bullet)}\notin\bigcup\Places^{(\nu)},\\
\nu= GE(A) & \biimplies & \bigcup
A^{(\bullet)}\in\bigcup\Places^{(\nu+1)}\setminus\bigcup\Places^{(\nu)},
\\
\nu> GE(A) & \biimplies & \bigcup
A^{(\bullet)}\in\bigcup\Places^{(\nu)}.
\end{array}
$$
Further elementary properties, whose proofs are left to the reader, are
stated in the next lemma.

\begin{mylemma}
Let
$\big(\Sigma_\mu\big)_{\mu\leqslant\xi},(\bullet),\TARGETS,\REDS,\POWNODES$
be a colored $\Places$-process and let $A \subseteq \Places$ be
a node. Then
\begin{itemize}
\item $A^{(\alpha)} = A^{(\bullet)}$, where $\alpha = GE(A)$;
\item if $q^{(\nu+1)} \supsetneq q^{(\nu)}$, for some $q \in A$
and some $\nu < \xi$, then $GE(A) > \nu$.
\end{itemize}
\end{mylemma}

Other important related definitions are the following.

\begin{mydef}\rm
A place $g$ is said to be a {\em local trash} for a node $A$ if
\begin{itemize}
\item $g\in\Targets{A}\setminus\REDS$, i.e., $g$ is a green
target of $A$;
\item there holds $GE(A)<GE(B)$, for every node $B$ such that $g\in B$.
\end{itemize}
\end{mydef}

\begin{mydef}
 A set $\mathcal{W}$ of places is said to be {\em closed} if
 \begin{itemize}
 \item all of its elements are green;
 \item every $\Pow$-node which intersects $\mathcal{W}$ has a local
 trash which belongs to $\mathcal{W}$. 
 \end{itemize}
\end{mydef}

\subsection{Minus-Surplus refinement}
In this section we recall some technical notions
to refine the original transitive partition.
This procedure stores some elements 
(the {\em Surplus} portion of a block) in order to trigger off a
construction which is supposed to ``pump''
elements inside fixed bocks. Conversely, the remaining
collection of elements (the {\em Minus} portion of a block)
will be used to copy the original formative
process.

We shall adopt the following notation. For a couple of
ordinals $\beta ',\beta ''$ we denote by 
$[\beta ',\beta '']$ the collection of ordinals
$\{\beta\mid \beta '\le\beta\le\beta ''\}$.

We say that a transitive partition $\Sigma$ is equipped of a Minus-Surplus 
partitioning if each block $q$ is partitioned into two sets, 
namely, $Surplus(q)$ and
$Minus(q)$. Consistently, we can extend this notation to a 
formative process $\big(\Sigma_\mu\big)_{\mu\leqslant\xi}$.
Given a node $\Gamma$, we indicate by $Minus(\Gamma ^{(\mu)})$ 
the collection
of sets 
$$\{Minus(q^{(\mu)})\mid q\in\Gamma\}.$$
Define now a Minus-Surplus partitioning for $\Sigma _0$, and assume that 
for each step $\mu$ of the process a refinement of the partition
$\{\Delta ^{(\mu)}(q)\}_{q\in\Sigma}$ is decided
in the following way:
for each $q\in\Sigma$ the set $\Delta ^{(\mu)}(q)$ 
is partitioned into two sets 
$\Delta ^{(\mu)}Minus(q)\subseteq\powAst (Minus(A_{\mu}^{(\mu)}))$
and $\Delta ^{(\mu)}Surplus(q)\subseteq(\powAst(A_{\mu}^{(\mu)})\setminus
\powAst (Minus(A_{\mu}^{(\mu)}))$.\\
Then define inductively $$Surplus(q^{(\mu+1)})=Surplus(q^{(\mu)})\cup 
\Delta ^{(\mu)}Surplus(q)$$ and
$$Minus(q^{(\mu+1)})=Minus(q^{(\mu)})\cup\Delta ^{(\mu)}Minus(q).$$
As far as $\xi$ limit are concerned, we put
$$Minus(q^{(\xi)})=\bigcup_{\mu <\xi}Minus(q^{(\mu)})$$
and, analogously,
$$Surplus(q^{(\xi)})=\bigcup_{\mu <\xi}Surplus(q^{(\mu)})$$
If $\Gamma$ is a subset of $\Sigma$, we denote by $Surplus(\Gamma)$ the 
set
$$\{q\mid q\in\Gamma\wedge Surplus(q)\neq\emptyset\}.$$
\begin{mydef}\label{}\rm 
Whenever a Surplus-Minus partition is defined for all blocks of a
transitive partition $\Sigma$, we say that $\Sigma$ is equipped of
a \emph{Minus-Surplus partitioning},
and we denote by $Surplus$-$Minus(\Sigma )$ the following
refinement of the original one:
$$\{Minus(q),Surplus(q)\mid q\in\Sigma\}.$$
It is rather obvious that $Surplus$-$Minus(\Sigma )\sqsubseteq\Sigma$.
\end{mydef}

\begin{myremark}\label{powdisj}\rm
Easy combinatorial arguments (see \cite[Lemma 3.1 5(b)]
{Cantone-Omodeo-Ursino02}) 
show that $\powast\_$ of Surplus
and Minus nodes are mutually disjoint.
\end {myremark}

The next definition says which structural properties 
a formative process has to fulfill in order to copy 
the history of a transitive partition.

\begin{mydef}\label{formimit}\rm 
Let
$(\left(\{q^{(\mu)}\}_{q\in\Places}\right)_{\mu\leqslant\xi},(\bullet),
\TARGETS,\REDS,\POWNODES)$
be a {colored $\Places$-process}. Besides, let
$\big(\{\widehat{q}^{[\alpha]}\}_{\widehat{q}\in\widehat{\Places}}\big)_
{\alpha\in [\alpha',\alpha '']}$
a formative processes equipped of a Minus-Surplus 
partitioning. Assume that
$q\rightarrow \widehat{q}$ is a bijection from $\Places$ to 
$\widehat{\Places}$,
$\beta ''\le\xi$, and $\gamma$ is an order preserving 
injection from
$[\beta ',\beta '']$ to $[\alpha',\alpha '']$.
Let ${\mathcal C}$ be a closed collection of green blocks,
and $q\rightarrow \widehat{q}$ be a bijection from $\Places$ to 
$\widehat{\Places}$.
We
say that $\big(\{\widehat{q}^{[\alpha]}\}_{\widehat{q}
\in\widehat{\Places}}\big) _{\alpha\in\gamma 
[[\beta ',\beta '']]}$ \emph{imitates}
the segment $[\beta ',\beta '']$ of the process 
$(\left(\{q^{(\mu)}\}_{q\in\Places}\right)_{\mu\leqslant\xi}$
if the following hold
for all $\beta$ in $[\beta ',\beta '']$:
\begin{itemize}
\item[(i)]$\vert q^{(\beta)}\vert=\vert Minus^{[\gamma(\beta)]}(\widehat{q})\vert$;
\item[(ii)] $\vert \Delta ^{(\beta)}(q)\vert=
\vert\Delta ^{[\gamma (\beta)]} Minus(\widehat{q})\vert$;
\item[(iii)]$\Delta ^{[\gamma (\beta)]}Surplus(\widehat{q})\neq\emptyset$ 
implies
$\beta=GE(A_{\beta})$, $q$ local trash for $A_{\beta}$ and
$q\in {\mathcal C}$;
\item[(iv)] If $\Gamma\in\POWNODES$ holds, then
$\powast{\widehat{\Gamma}^{[\gamma (GE(\Gamma))]} }
\subseteq\bigcup\widehat{\Sigma}^{[\gamma (GE(\Gamma)+1)]}$;
\item[(v)]For all $\beta\ne GE(\Gamma)$  
$\bigcup\Gamma^{(\beta)}\in\Delta ^{(\beta)}(q)$ iff
$\bigcup Minus\widehat{\Gamma}^{[\gamma (\beta)]}\in\Delta ^{[\gamma (\beta)]}(
\widehat{q})$;
\item[(vi)] If $\beta =GE(\Gamma)$ then $\bigcup\Gamma^{(\beta)}\in 
\Delta ^{(\beta)}(q)$ iff
$\bigcup\Gamma^{[\gamma (\beta)]}\in\Delta ^{[\gamma (\beta)]}(\widehat{q})$;
\item[(vii)] For all $q\in\REDS$ $\widehat{q}^{[\gamma (\beta)]}=
Minus^{[\gamma(\beta)]}(\widehat{q})$;
\item[(viii)] For all ordinals $\beta$ $\{q\mid 
\widehat{q}\in Surplus(\widehat{\Sigma})^{[\gamma (\beta)]}\}
\subseteq{\mathcal C}$;
\item[(ix)] $\vert\powAst (Minus (\widehat{\Gamma}^{[\gamma(k)]})) 
      \setminus \bigcup _{q\in\Sigma}q^{[\gamma(k)]}\vert=
\vert\powAst (\Gamma) ^{(k)}\setminus \bigcup _{q\in\Sigma}q^{(k)}\vert$;
\item[(x)] $\vert\powAst (Minus(\widehat{\Gamma}^{[\gamma(k-1)]})) \cap  q^{ 
[\gamma(k)]}\vert  
=\vert \powAst (\Gamma^{(k-1)})\cap q^{(k)}\vert$. 
\end{itemize}
\end{mydef}

\smallskip
\begin{myremark}\rm
\label{rem1}
We make some simple observations.
\begin{itemize}
\item $\powAst (\Gamma^{(k-1)})\cap q^{(k)}=\powAst (\Gamma^{(k)})\cap q^{(k)}$. Hence,
whenever $\gamma (k)$ is the successor of $\gamma (k-1)$, (x) can be rephrased as
$$\vert\powAst (Minus(\widehat{\Gamma}^{[\gamma(k)]})) \cap  q^{ 
[\gamma(k)]}\vert =\vert \powAst (\Gamma^{(k)})\cap q^{(k)}\vert.$$
\item Naturally, (ix) belongs to the structural properties that a 
formative process has to fulfill in order to
simulate another one, although it can be obtained from (i) and (x). 
\item Assume that (viii) holds at the beginning of the process.
Then (iii) entails (viii), therefore, whenever
one has to prove inductively
the previous properties, it suffices to show that (viii) holds only
in the starting step.
The same argument holds for (x).
Indeed, it can be obtained from (ii), (iii) and (x) of the preceding step.
\end{itemize}
\end {myremark}


The following requirements set are to be satisfied by
the initial conditions of a transitive partition 
in order to
play the role of starting point of an imitation process
(as it is easily seen, they are purely combinatorial).

\begin{mydef}\label{weakimit}\rm 
Let
$(\left(\{q^{(\mu)}\}_{q\in\Places}\right)_{\mu\leqslant\xi},
(\bullet),\TARGETS,\REDS,\POWNODES)$
be a {colored $\Places$-process},
$(\widehat{\Sigma},\widehat{\mathcal{G}})$ be a {$\widehat{\Sigma}$-board}
equipped with a Minus-Surplus 
partitioning,
$q\rightarrow \widehat{q}$ be a bijection from $\Places$ 
to $\widehat{\Places}$, and 
${\mathcal C}$ be a closed collection of green blocks.
Assume $k'<\xi$,
such that (i), (vii), (viii) and (x) of Def.\ref{formimit} hold in the version
$\widehat{\Sigma}_{\gamma (k')}=
\widehat{\Sigma}$.
We say that $\widehat{\Sigma}$ \emph{weakly imitates} 
$\Sigma$ \emph{upwards}, provided that
the following conditions are satisfied:
\begin{itemize}
\item [(a)] for all $\Gamma\subseteq\Sigma$ and $q\in\Sigma$, 
$$\bigcup Minus(\widehat{\Gamma})\in\powAst (Minus(\widehat{\Gamma})) 
\setminus \bigcup _{q\in\Sigma}\widehat{q}\quad {\rm iff} 
\quad\bigcup\Gamma^{(k')}\in
\powAst (\Gamma^{(k')})\setminus \bigcup _{q\in\Sigma}q^{(k')};$$
\item [(b)] $q\in\Gamma\wedge Surplus(q)\neq\emptyset\wedge GE(\Gamma)\ge k'$
implies $\bigcup\widehat{\Gamma}\in\powAst 
(\widehat{\Gamma}) 
\setminus \bigcup _{q\in\Sigma}\widehat{q}$;
\item [(c)] if $GE(\Gamma)<k'$, then $\bigcup\Gamma^{(k')}\in q^{(k')}$ iff
$\bigcup \widehat{\Gamma}\in \widehat{q}$ and $\Gamma\in\POWNODES$
implies $\powAst (\widehat{\Gamma})\subseteq\bigcup\widehat{\Sigma}$.
\end{itemize}
\end{mydef}

\section{Two structural results concerning Minus-Surplus partition}


The following Lemma relates Definition \ref{weakimit} 
with the notion of imitating a formative process.

\begin{mylemma}\label{pasting}
Let
$(\left(\{q^{(\mu)}\}_{q\in\Places}\right)_{\mu\leqslant\xi},
(\bullet),\TARGETS,\REDS,\POWNODES)$
be a {colored $\Places$-process},
$(\widehat{\Sigma},\widehat{\mathcal{G}})$ be a {$\widehat{\Sigma}$-board},
the latter equipped of a Minus-Surplus 
partitioning,
$q\rightarrow \widehat{q}$ be a bijection from $\Places$ 
to $\widehat{\Places}$, and
${\mathcal C}$ be a closed collection of green blocks.
Assume that $k'\le\xi$, and that  
$\widehat{\Sigma}$ weakly imitates upward $\Sigma_{k'}$.
Define $\widehat{\Sigma}=\widehat{\Sigma}_{\gamma(k')}$ and, for all 
$q\in\widehat{\Places}$,
$\widehat{q}=\widehat{q}^{[\gamma(k')]}$. 
Then for all ordinals 
$k''$ such that $k''\le\xi$ and $\vert [k',k'']\vert <\omega$
it can be constructed a formative process 
$\big(\{\widehat{q}^{[\alpha]}\}_{\widehat{q}
\in\widehat{\Places}}\big)_{\gamma(k')\leqslant\mu\leqslant\gamma(k'')}$
which imitates the
segment $[k',k'']$ of the process $(\left(\{q^{(\mu)}\}_{q\in\Places}\right)_{\mu\leqslant\xi},
(\bullet),\TARGETS,\REDS,\POWNODES$.
\end{mylemma} 

\begin{proof}
We construct a formative process by induction 
satisfying the requested properties (i)-(x).

Concerning the base case $\mu =\gamma (k')$, 
(i),(vii),(viii)(x) hold by hypothesis, and
(ix) holds by Remark~\ref{rem1}, since (i) and (x) hold.
Assume $k'\ne GE(A_{k'})$.
Using (ix) and hypothesis (a) we can define a partition 
$\bigcup_{q\in\Sigma}(\Delta
^{[\gamma(k')]}(q)$ of 
$$\powAst (Minus
^{[\gamma(k')]}(\widehat{A_{k'}}))\setminus\bigcup _{q\in\Sigma}q^{[\gamma(k')]}$$
such that (ii) and (v) hold, as well.
If $k'=GE(A_{k'})$ and $Surplus(\widehat{q}^{[\gamma(k')]})\neq\emptyset
$ for some $q\in A_{k'}$ 
(otherwise we proceed as before, and condition
(vi) is automatically satisfied),  then, 
using (b), interchanging $\bigcup Minus(\widehat{A_{k'}}^{[\gamma(k')]})$ with 
$\bigcup A^{[\gamma(k')]}$, (vi) is satisfied.

If $A_{k'}\in\POWNODES$ and $\widehat{A_{k'}}=Minus(\widehat{A_{k'}})$, proceed as
before (in this case (iv) holds by a straight checking 
of cardinality starting from (ix)).
Otherwise, since (viii) holds, 
there must exist a local trash 
$q\in {\mathcal C}$ for $A_{k'}$.
Then, construct the partition as before, except
for $\Delta^{[\gamma(k')]}Surplus(\widehat{q})$, in which  we put 
the whole remainder 
$$(\powAst (\widehat{A_{k'}}^{[\gamma(k')]}) 
\setminus \bigcup _{q\in\Sigma}q^{[\gamma(k')]})\setminus
\bigcup_{q\in\Sigma}\Delta^{[\gamma(k')]} Minus(\widehat{q}),$$
so satisfying
(iii) and (iv).

Now, assume all the inductive hypotheses for $\gamma(k)$. Our aim 
is to demonstrate
the case $\gamma(k+1)$.
By Remark~\ref{rem1}, provided that (iii)[$\gamma(k+1)$] is proven, 
(viii) automatically holds.
Plainly (i)[$\gamma(k)$] and (ii)[$\gamma(k)$] entail
(x)[$\gamma(k+1)$] and (i)[$\gamma(k+1)$]. The latter
in turns implies the following for all $\Gamma\subseteq\Sigma$
\begin{equation}
\label{pow2}
\vert\powAst (Minus ^{[\gamma(k+1)]}\widehat{\Gamma}))\vert=\vert
\powAst (\Gamma^{(k+1)})\vert.
\end{equation}

In order to show (ix) we observe that,
since
\begin{align*}
&\powAst (Minus(\widehat{\Gamma}^{[\gamma(k+1)]})\setminus
\bigcup _{q\in\Sigma}q^{[\gamma(k+1)]}\\
&= (\powAst (Minus(\widehat{\Gamma}^{[\gamma(k+1)]})
\setminus\powAst (Minus(\Gamma^{[\gamma(k)]}))
\setminus 
\bigcup _{q\in\Sigma}q^{[\gamma(k+1)]})\cup \\
&\quad \cup(\powAst (Minus(\widehat{\Gamma}^{[\gamma(k)]}))
\setminus 
\bigcup _{q\in\Sigma}q^{[\gamma(k+1)]}),
\end{align*}
it follows that
\begin{align*}
&\powAst (Minus(\widehat{\Gamma}^{[\gamma(k+1)]}))\setminus
\powAst (Minus(\widehat{\Gamma}^{[\gamma(k)]}))
\setminus
\bigcup _{q\in\Sigma}q^{[\gamma(k+1)]}\\
&=\powAst (Minus(\widehat{\Gamma}^{[\gamma(k+1)]}))
\setminus\powAst (Minus(\Gamma^{[\gamma(k)]})).
\end{align*}
Therefore,
\begin{align*}
&\powAst (Minus(\widehat{\Gamma}^{[\gamma(k+1)]}))\setminus
\bigcup _{q\in\Sigma}q^{[\gamma(k+1)]}\\
&=\powAst (Minus(\widehat{\Gamma}^{[\gamma(k+1)]}))
\setminus\powAst (Minus(\widehat{\Gamma}^{[\gamma(k)]}))
\cup\powAst (Minus(\widehat{\Gamma}^{[\gamma(k)]}))
\setminus 
\bigcup _{q\in\Sigma}q^{[\gamma(k+1)]}.
\end{align*}
Reasoning in the same way, we obtain 
\begin{align*}
&\powAst (\Gamma)^{(k+1)}\setminus\bigcup_{q\in\Sigma}q^{(k+1)}\\
&= \powAst (\Gamma)^{(k+1)}\setminus\powAst (\Gamma)^{(k)}\cup\powAst (\Gamma)^{(k)}
\setminus\bigcup _{q\in\Sigma}q^{(k+1)}.
\end{align*}
By the induction hypothesis (i)[$\gamma(k)$] we have 
$\vert\powAst (Minus ^{[\gamma(k)]}\widehat{\Gamma}))\vert=\vert
\powAst (\Gamma^{(k)})\vert$, and by 
equation~\eqref{pow2},
$$\vert\powAst (Minus ^{[\gamma(k+1)]}\widehat{\Gamma}))\vert=\vert
\powAst (\Gamma^{(k+1)})\vert,$$
which in turns implies
$$\vert\powAst (Minus(\widehat{\Gamma}^{[\gamma(k+1)]}))\setminus\powAst 
(Minus(\widehat{\Gamma}^{[\gamma(k)]}))\vert=\vert\powAst (\Gamma)^{(k+1)}\setminus
\powAst (\Gamma)^{(k)}\vert.$$

Hence we are left to prove the equality
\begin{equation}
\label{pow3}
\vert\powAst (Minus(\widehat{\Gamma}^{[\gamma(k)]}))
\setminus 
\bigcup _{q\in\Sigma}q^{[\gamma(k+1)]}\vert=\vert\powAst (\Gamma)^{(k)}
\setminus\bigcup _{q\in\Sigma}q^{(k+1)}\vert.
\end{equation}
Observe that 
$$\powAst (Minus(\widehat{\Gamma}^{[\gamma(k)]}))
\setminus 
\bigcup _{q\in\Sigma}q^{[\gamma(k+1)]}=\powAst (Minus(\widehat{\Gamma}^{[\gamma(k)]}))
\setminus\bigcup _{q\in\Sigma}
q^{[\gamma(k)]}\setminus\bigcup _{q\in\Sigma}\Delta^{[\gamma (k)]}Minus(q).$$
If $\Gamma\neq A_{k}$, by the disjointness of $\powAst$ we get
$$\powAst (Minus(\widehat{\Gamma}^{[\gamma(k)]}))
\setminus\bigcup _{q\in\Sigma}
q^{[\gamma(k)]}\setminus\bigcup _{q\in\Sigma}\Delta^{[\gamma (k)]}Minus(q)=
\powAst (Minus(\widehat{\Gamma}^{[\gamma(k)]}))\setminus\bigcup _{q\in\Sigma}
q^{[\gamma(k)]}.$$
Plainly, the same is true in the $\_ ^{()}$ version, thus \eqref{pow3}
holds for $\gamma(k)$, by virtue of (ix). 
Otherwise, since $\bigcup _{q\in\Sigma}\Delta^{[\gamma (k)]}Minus(q)$
is a partition of a subset extract from
$$\powAst (Minus(\widehat{\Gamma}^{[\gamma(k)]}))
\setminus 
\bigcup _{q\in\Sigma}q^{[\gamma(k)]},$$
we have that
\begin{align*}
&\vert\powAst (Minus(\widehat{\Gamma}^{[\gamma(k)]}))
\setminus\bigcup _{q\in\Sigma}
q^{[\gamma(k)]}\setminus\bigcup _{q\in\Sigma}\Delta^{[\gamma (k)]}
Minus(q)\vert\\
&=\vert\powAst (Minus(\widehat{\Gamma}^{[\gamma(k)]}))
\setminus\bigcup _{q\in\Sigma}
q^{[\gamma(k)]}\vert -\sum _{q\in\Sigma}\vert\Delta^{[\gamma (k)]}
Minus(q)\vert.
\end{align*}
Again, the same holds in the $\_ ^{()}$ version, and 
\eqref{pow3} is reached by (i)[$\gamma(k)$] and (ii)[$\gamma(k)$].
This concludes the proof of (ix)[$\gamma(k+1)$].

Concerning (vii)[$\gamma(k+1)$], observe that 
$q^{[\gamma(k+1)]}=q^{[\gamma(k)]}\cup\Delta^{[\gamma (k)]}(q)$.
By the induction hypothesis (vii)[$\gamma(k)$],
$$q^{[\gamma(k)]}=Minus(q^{[\gamma(k)]}).$$ 
On the other side,
since (iii)[$\gamma(k)$] holds and ${\mathcal C}$
is composed of green places only, 
$$\Delta^{[\gamma (k)]}(q)=\Delta^{[\gamma (k)]}Minus(q),$$
which implies (vii)[$\gamma(k+1)$].

Regarding (ii)[$\gamma(k+1)$]-(vi)[$\gamma(k+1)$], the argument goes like 
in the base case.
\end{proof}

\begin{mylemma}\label{upward}
Let 
$(\left(\{q^{(\mu)}\}_{q\in\Places}\right)_{\mu\leqslant\xi},
(\bullet),\TARGETS,\REDS,\POWNODES)$
be a {colored $\Places$-process}. Moreover, let
$\big(\{\widehat{q}^{[\alpha]}\}_{\widehat{q}\in
\widehat{\Places}}\big)_{\alpha\leqslant\xi '}$ be another formative process,
equipped of a Minus-Surplus partitioning.
Assume that, for some $k'\le\xi$ and $m\le\xi '$, 
\begin{itemize}
\item $\widehat{\Sigma}_{m}$ {weakly imitates}
$\Sigma_{k'}$ {upwards};
\item the process
$\big(\{\widehat{q}^{[\alpha]}\}_{\widehat{q}\in\widehat{\Places}}
\big)_{\alpha\in\gamma [k',\xi]}$ {imitates}
$\big(\{q^{(\mu)}\}_{q\in\Places}\big)_{k'\le\mu\leqslant\xi}$, where $\gamma$ is an
injective map from $[k',\xi]$ to $[m,\xi ']$;
\item $\widehat{\Sigma}_{\xi '}$ has the same targets of $\Sigma_{\xi}$;
\item for all $\mu >m\wedge\mu\notin\gamma [k',\xi]$ the following holds:
$\Delta^{[\mu]}(\widehat{q})\subseteq\Delta^{[\mu]}Surplus(\widehat{q})$;
\item 
if $\beta$ is the greatest ordinal
such that
$\beta\in\gamma[k',\xi]\wedge\beta\le\mu$, if 
$q$ is a local trash of ${A}_{\mu}$, and if $GE(A_{\mu})>\gamma^{-1}(\beta)$,
then $\bigcup\widehat{A_\mu}^{[\xi ']}\notin
\Delta^{[\mu]}Surplus(\widehat{q})$.
\end{itemize}
Then $\widehat{\Sigma}_{\xi '}$ {imitates} $\Sigma_{\xi}$ {upwards}.
\end{mylemma}

\begin{proof}
We prove that the resulting partition $\widehat{\Sigma}_{\xi '}$
fulfills the 
conditions:
\begin{itemize}
\item[(0)]  $q^{(\xi)}\ni\in\powast{\Gamma}^{(\xi)}$
holds if and
only if $\widehat{q}^{[\xi ']}\ni\in\powast{\widehat{\Gamma}^{[\xi ']}}$;
\item[(1)] $\bigcup\widehat{\Gamma}^{[\xi ']}\in\widehat{q}^{[\xi ']}$ 
if and only if $\bigcup \Gamma^{(\xi)}\in q^{(\xi)}$;
\item[(2)] if $\Gamma\in\POWNODES$ holds, then
$\powast{\widehat{\Gamma}^{[\xi ']}}\subseteq\bigcup\widehat{\Places}^{[\xi ']}
$;
\item[(3$'$)] if
$q\in\REDS$, then $\card{\widehat{q}^{[\xi ']}}=\card{q^{(\xi)}}$.
\end{itemize}

Along the verification of properties (0)-(3$'$) we refer to (i)-(x) 
of Def.~\ref{formimit}.
\begin{itemize}
\item[(0)] By the fact that the two partitions have the same targets;

\item[(1)] 
In case $\bigcup\widehat{\Gamma}^{[\xi ']}\in\widehat{q}^{[\xi ']}$, 
assuming that it
is distributed strictly before $m$, then $GE(\Gamma)<k'$. Indeed, 
if not so, by (vi) Def.~\ref{formimit}, since $\bigcup\Gamma^{(GE(\Gamma))}\in 
\Delta ^{(GE(\Gamma))}(q)$,
$$\bigcup\widehat{\Gamma}^{[\xi ']}\in\widehat{q}^{[\xi ']}=
\bigcup\widehat{\Gamma}^{[\gamma (GE(\Gamma))]}
\in\Delta^{[\gamma (GE(\Gamma))]}\widehat{q},$$
which is impossible, due to the fact
that  $\bigcup\widehat{\Gamma}^{[\xi ']}\in\widehat{q}^{[\xi ']}$
is already in $\widehat{q}^{[\gamma (GE(\Gamma))]}$, and 
$\Delta^{[\gamma (GE(\Gamma))]}\widehat{q}$,
by definition, is made of elements of 
$\powAst (\widehat{\Gamma}^{[\gamma(GE(\Gamma)]})
\setminus 
\bigcup _{q\in\Sigma}q^{[\gamma(GE(\Gamma))]}$.
Then, using the fact that 
$\widehat{\Sigma}_{m}$ 
weakly simulates  $\Sigma_{k'}$, the result follows.
Concerning the right implication, we are left to prove the case when
$\bigcup\widehat{\Gamma}^{[\xi ']}$
is distributed after or in $m$.
Let $j$ be such an index. By hypothesis,
$j$ cannot be outside $\gamma [k',\xi]$, and so
$j=\gamma (k)$ for some $k$. 
We show that $k=GE(\Gamma)$. By contradiction, let us assume $k>GE(\Gamma)$. 
Then, by (vi) Def.\ref{formimit},
$$\bigcup\widehat{\Gamma}^{[\gamma (GE(\Gamma))]}
\in\Delta^{[\gamma (GE(\Gamma))]}\widehat{q}.$$
Observe that, after $\gamma (GE(\Gamma))$, $\widehat{\Gamma}$ cannot change inside
the range of $\gamma$, on
account of (ii) and (iii) of Def.\ref{formimit}. It
it cannot change for 
an index $j$ outside, since
$GE(\Gamma)$ is greater than the greatest ordinal $\beta$ such that
$\beta\in\gamma[k',\xi]\wedge\beta\le j$.
On the other hand, 
$k$ cannot be strictly less than $GE(\Gamma)$, since in this case the same argument
used for $\bigcup\widehat{\Gamma}^{[\xi ']}$ distributed before $m$ and $GE(\Gamma)\ge k'$
applies.
Therefore $k=GE(\Gamma)$, and we are done.
We now show the left implication in the case $GE(\Gamma)< k'$.
The hypothesis implies  that
$\bigcup\widehat{\Gamma}^{[m]}\in\widehat{q}^{[m]}$. Reasoning as before,
we conclude that $\widehat{\Gamma}$ cannot change along the process after $m$.
Finally, assuming $GE(\Gamma)\ge k'$, by (vi) \eqref{formimit} there holds
$$\bigcup\widehat{\Gamma}^{[\gamma (GE(\Gamma))]}\in
\Delta^{[\gamma (GE(\Gamma))]}\widehat{q}.$$
Again $\widehat{\Gamma}$ cannot change in the sequel of the
process, either along
the imitated process, or outside.

\item[(2)]Follows plainly from (iv) \eqref{formimit}. Indeed, 
$\Gamma\in\POWNODES$, therefore
$$\powast{\widehat{\Gamma}^{[\gamma (GE(\Gamma))]}}
\subseteq\bigcup\widehat{\Places}^{[\gamma (GE(\Gamma))]}.$$ 
As observed in the previous point, 
after $[\gamma (GE(\Gamma))]$, $\widehat{\Gamma}$ cannot change either
along the imitating process, by (ii) and (iii) \eqref{formimit}, 
or outside, by hypothesis.
Thus $\powast{\widehat{\Gamma}^{[\xi ']}}\subseteq\bigcup\widehat{\Places}^{[\xi ']}$.

\item[(3$'$)]The red places cannot belong to  
${\mathcal C}$. Hence, by the property (viii),
they cannot have Surplus part, which in turns implies that
$Minus(\widehat{q}^{[\xi ']})=\widehat{q}^{[\xi ']}$.
This, combined with
$\vert Minus^{[\gamma(\xi)]}(q)\vert=\vert q^{(\xi)}\vert$, due to 
(i) \eqref{formimit}, 
leads to the thesis.
\end{itemize}
\end{proof}
The following theorem summarizes the previous results and
shows which properties two formative processes have to share in order to model the same literals.
The proof is a straight application of Corollary \ref{MsatisfiesII}
\begin{mytheorem}\label{upwards}
Let 
$(\left(\{q^{(\mu)}\}_{q\in\Places}\right)_{\mu\leqslant\xi},
(\bullet),\TARGETS,\REDS,\POWNODES)$
be a {colored $\Places$-process}. Moreover, let
$\big(\{\widehat{q}^{[\alpha]}\}_{\widehat{q}\in
\widehat{\Places}}\big)_{\alpha\leqslant\xi '}$ be another formative process,
equipped of a Minus-Surplus partitioning.
Assume that, for some $k'\le\xi$ and $m\le\xi '$, 
\begin{itemize}
\item $\widehat{\Sigma}_{m}$ {weakly imitates}
$\Sigma_{k'}$ {upwards};
\item the process
$\big(\{\widehat{q}^{[\alpha]}\}_{\widehat{q}\in\widehat{\Places}}
\big)_{\alpha\in\gamma [k',\xi]}$ {imitates}
$\big(\{q^{(\mu)}\}_{q\in\Places}\big)_{k'\le\mu\leqslant\xi}$, where $\gamma$ is an
injective map from $[k',\xi]$ to $[m,\xi ']$;
\item $\widehat{\Sigma}_{\xi '}$ has the same targets of $\Sigma_{\xi}$;
\item for all $\mu >m\wedge\mu\notin\gamma [k',\xi]$ the following holds:
$\Delta^{[\mu]}(\widehat{q})\subseteq\Delta^{[\mu]}Surplus(\widehat{q})$;
\item 
if $\beta$ is the greatest ordinal
such that
$\beta\in\gamma[k',\xi]\wedge\beta\le\mu$, if 
$q$ is a local trash of ${A}_{\mu}$, and if $GE(A_{\mu})>\gamma^{-1}(\beta)$,
then $\bigcup\widehat{A_\mu}^{[\xi ']}\notin
\Delta^{[\mu]}Surplus(\widehat{q})$.
\end{itemize}
Consider a formula $\Phi\in {\rm MLSSPF}$,
a set-valued assignment ${\mathcal M}\in\{\:\mbox{\rm sets}\:\}^{{\mathcal X}_{\Phi}}$ 
defined on the collection ${\mathcal X}_{\Phi}$
of variables in $\Phi$
assuming that $(\left(\{q^{(\mu)}\}_{q\in\Places}\right)_{\mu\leqslant\xi},
(\bullet),\TARGETS,\REDS,\POWNODES)$
is a {colored $\Places$-process} for the \emph{$\Sigma_{{\mathcal X}_{\Phi}}$-board}

then, letting ${\mathcal M}'(v)=\bigcup[\Im_{{\mathcal M}}(\widehat{v})]$, 
for every literal in $\Phi$, the following conditions are
fulfilled:
\begin{des}
\item if the literal is satisfied by $\mathcal M$, then it is satisfied
by ${\mathcal M}'$ too;
\item if the literal is satisfied by ${\mathcal M}'$, and does not
involve $\Pow$~ or the construct $\{\anonymous,\dots,\anonymous\}$,
then it is satisfied by ${\mathcal M}$ too.
\end{des}
\end{mytheorem}

\begin{myremark}\rm
\label{salient}
The same result holds even in more relaxed conditions, revealing its
strength when we are looking for small models.
Namely, when we prune
the process instead of prolongate it.
In fact, the previous theorem holds, with an identical proof, 
provided that the domain of $\gamma$ 
contains the following two collections of
{\em salient} ordinals:
$$M_{arrow}=\{ \mu\mid k'\le\mu<\xi\wedge
\exists q\in \Places
q^{(\mu)} \cap \powast{A^{(\mu)}_{\mu}} =\emptyset\wedge
\Delta^{(\mu)}(q)\neq\emptyset\}$$
and
$$M_{GE}=\{\mu\mid k'\le\mu<\xi
\wedge
\bigcup A_{\mu}^{(\mu)} =\bigcup A_{\mu}^{(\bullet)}
\in\bigcup P^{(\bullet)}\}.$$
\end {myremark}

\section{Using the Two Structural Lemmas into Set Computable Examples}

Assume that ${\mathcal M}$ is a finite set assignment to the 
variables of an assigned formula $\Phi$ of MLSSPF, which
contains literals of the type $\neg Finite(x)$. 
Obviously, ${\mathcal M}$ cannot be a model for $\Phi$, although it
could happen that it satisfies every other literal, except those 
of that kind. The question is: in this situation could 
${\mathcal M}$ witness the satisfiability of $\Phi$? 
The answer is positive, as we will show, and
the core argument for proving that
lies inside a possible history of ${\mathcal M}$.
Indeed, given a formative process for the Venn partition 
$\Sigma$ inherited from  ${\mathcal M}$,
if we can find an ``engine" capable to pump elements inside at least one Venn region 
for each variable $x$, such that $\neg Finite(x)$
lies in $\Phi$ without affecting the satisfiability of other literals, we reach
the desired conclusion.

We will be more precise on the exact meaning
of ``engine'', and how profitably the results of the previous 
sections can be used in order to preserve the satisfiability of the other literals.
even though the
size of the assignment of some variables is infinitely increased. 

\begin{mydef}\rm
\label{prepareForPumpingPaths} In a $\Places$-board
$\mathcal{G}$, a \emph{path} is an ordered vertex list
$W_1,\dots,W_k$, in which places and nodes are so alternate that
$W_i,W_{i+1}$ is an edge of $\mathcal{G}$, for $i=1,\dots,k-1$. A
path is said to be \emph{simple} if neither places nor nodes occur
twice (i.e., $W_i\neq W_j$ when $0<i<j\leqslant k$ and
$i\equiv j\;(\mbox{{mod }}2)$ ).
\end{mydef}

\begin{mydef}\rm
\label{pumpingPaths} In a colored $\Places$-board
$\mathcal{G}=(\TARGETS,\REDS,\POWNODES)$, a path
$$\mathcal{C} \equiv
C_0,q_0,C_1,\dots,q_n,C_{n+1} $$
where the piece $C_0,q_0,C_1,\dots,q_n,$ is simple and
$n\geqslant 0$, \emph{devoid} of red places, and such
that $C_{n+1} = C_{0}$,
is said to be a
\emph{simple pumping cycle}.
\end{mydef}

Given a path $\mathcal{D}$ in a $\Places$-board $\mathcal{G}$, we
denote by $\PLACES{\mathcal{D}}$ and $\NODES{\mathcal{D}}$ the
collections of places and nodes occurring in $\mathcal{D}$, respectively.
Moreover, given a node $B$ in $\mathcal{G}$, we denote
by $\PN{B}$ the collection of all nodes which have nonnull
intersection with $B$.

The following is to be regarded as the engine which increases the size
of some places without affecting the validity of the formula.

\begin{mydef}\rm
\label{pumpchain}
\rm Let $\mathcal{C}$ be a simple pumping cycle
relative to a given colored $\Places$-process
$\big(\Sigma_\mu\big)_{\mu\leqslant\ell},[\bullet],\TARGETS,\REDS,\POWNODES$,
with $\ell$ finite. Then $<q_0,i_0,\mathcal{C}>$ is called a
\emph{simple pumping event} whenever we have
\begin{itemize}
     \item[(i)] $q_0^{[i_{0}]} \setminus \bigcup\bigcup
     \Places^{[i_{0}]} \neq \emptyset$, $q_0\in
     \PLACES{\mathcal{C}}$;

     \item[(ii)]  $GE(\PN{\PLACES{\mathcal{C}}})\ge i_0$;

     \item[(iii)] $\powAst(B^{[i_{0}]}) \neq \emptyset$
     (i.e., $\emptyset\notin B^{[i_{0}]}$), for $B
     \in\NODES{\mathcal{C}}$.  
\end{itemize}
\end{mydef}

If $\Sigma$ is a particular Venn partition $\Sigma_{{\mathcal X}_{\Phi}}$,
the variables that contain the places involved in the pumping cycle can be considered
\emph{potential infinite variables}.

\begin{mytheorem}\label{pumping1} 
Assume that ${\mathcal M}$ is a finite transitive set 
assignment to the variables ${\mathcal X}_{\Phi}$ of an assigned formula 
$\Phi$ of {\rm MLSSPF}, that
satisfies every other literals except those
of the type $\neg Finite(x)$. Consider the transitive
\emph{$\Sigma_{{\mathcal X}_{\Phi}}$-board} ${\mathcal G}=(\TARGETS,\REDS,\POWNODES)$,
and an associated colored $\Places$-process
$\big(\Sigma_\mu\big)_{\mu\leqslant\ell},(\bullet),\TARGETS,\REDS,\POWNODES$,
with $\ell$ finite. 
Then there exists a model for $\Phi$,
provided there is a {simple pumping event}
$<q,i_0,\mathcal{C}>$ such that $\PLACES{\mathcal{C}}$ is 
contained in a closed set $\overline {\mathcal{C}}$
satisfying the statement:
$$\mbox{ For each variable } x \mbox{ such that }\neg Finite(x)\in\Phi,
\Im_{{\mathcal M}}(x)\cap\PLACES{\mathcal{C}}\mbox{ is not empty. }$$ 
\end{mytheorem}

\proof
Let $<q_0,i_0,\mathcal{C}>$ be our {simple pumping event}, where 
$\mathcal{C}$ is equal to $$\{C_0,q_0\dots q_{n}, C_{n+1}\}.$$ 
We build a new formative process 
$\big(\widehat{\Sigma}_\mu\big)_{\mu\leqslant\ell},[\bullet],\TARGETS $,
using the original one as an oracle.
In the meanwhile, a Minus-Surplus refinement is done.
We first define the sequence of the nodes to be used in this new process.
Denote by  $\seq {\ell}=\{A_0 \dots A_{\ell}\}$ the sequence of nodes
used along the given process $\big(\Sigma_\mu\big)_{\mu\leqslant\ell},
(\bullet),\TARGETS,\REDS,\POWNODES$.
The following sequence serves to our scope:
$$A_1,\dots A_{i_0 -1}\underbrace{C_1\dots C_{n+1}}_{\aleph_{0} -times}, 
A_{\gamma(i_0)}\dots A_{\gamma(\ell)},$$
where, for all $j$, $A_{\gamma(j)}=A_j$ and the cycle 
$\mathcal{C}$ are repeated $\aleph_{0} -times$.

In order to define a formative process, we just need to exhibit the way to
distribute all the elements produced at each stage.
Our strategy consists to follow the old formative process up to the
stage $i_0 -1 =\gamma (i_0-1)$, 
setting $\left(\{\widehat{q}^{[j]}\}_{q\in\Places}\right)_{j\leqslant i_0-1}$
$\left(\{q^{(j)}\}_{q\in\Places}\right)_{j\leqslant i_0-1}$.
Along this segment, 
we define $\gamma$ as the identity map; then, we ``pump'' the cycle in order to create
new elements and distribute them. 
This procedure by transfinite induction
increases the cardinality of the blocks inside the cycle, preserving the cardinality
of all the blocks not involved in the pumping procedure.
In order to do that, we 
distinguish the elements reserved for the pumping procedure (Surplus portion) from those
used for mimicking the old process (Minus portion). The
Minus-Surplus refinement that we are about
to define will serve such a scope.

Without loss of generality, 
we assume that at each step the cycle can distribute at least 
three new elements (otherwise, we can pump the cycle
to give at least two elements to every block involved in the cycle).
By  Definition of {simple pumping event}, $q^{(i_{0})}\setminus\bigcup\bigcup
\Places^{(i_{0})} \neq \emptyset$, which means that in $q^{(i_{0})}$ there are unused elements.
Let $t_0$ be one of these, and 
define the partitions Surplus and Minus as follows:
\begin{itemize}
\item For all $q\neq q _0$ put\\
      $Surplus ^{[\gamma (i_0-1)+1]}(\widehat{q}) =\emptyset$ and 
      $Minus^{[\gamma (i_0-1)+1]}(\widehat{q})=q^{(i_0)}$;
\item For $q _0$ put\\
      $Surplus ^{[\gamma (i_0-1)+1]} (\widehat{q} _0)=\{t_0\}$
      $ Minus^{[\gamma (i_0-1)+1 ]}(\widehat{q} _0)=
      q_0^{(i_0)}\setminus\{ t_0 \} $;
\end{itemize}
Since every block involved in the cycle has at least two elements, 
the set 
$$\powAst \Big( \big\{ Surplus^{[\gamma (i_0-1)+1]}(\widehat{q}_0)\big\} 
\cup \widehat{C_1}^{[\gamma (i_0-1)+1]}\Big)\setminus
\Big\{  \bigcup \widehat{C_1}^{[\gamma (i_0-1)+1]}\Big\}$$
is not empty.
Moreover, by Lemma \ref{unused}, it is made of unused elements only.
Thus,
\begin{align*}
&\powAst ( \{ Surplus^{[\gamma (i_0-1)+1]}(\widehat{q}_0)\} 
\cup\widehat{C_1} ^{[\gamma (i_0-1)+1]})\setminus
\{  \bigcup\widehat{C_1} ^{[\gamma (i_0-1)+1]}\}\\
&=(\powAst ( \{ Surplus^{[\gamma (i_0-1)+1]}(\widehat{q}_0)\} 
\cup\widehat{C_1} ^{[\gamma (i_0-1)+1]})\setminus
\{  \bigcup\widehat{C_1} ^{[\gamma (i_0-1)+1]}\})\setminus
\bigcup _{q\in\Sigma}\widehat{q}^{[\gamma(i_0-1)+1]},
\end{align*}
so that the position
$$\Delta ^{[\gamma (i_0-1)+1]}(Surplus(\widehat{q}_1))=
\powAst ( \{ Surplus^{[\gamma (i_0-1)+1]}(\widehat{q}_0)\} 
\cup \widehat{C_1}^{[\gamma (i_0-1)+1]})\setminus
\{  \bigcup \widehat{C_1}^{[\gamma (i_0-1)+1]}\}$$
makes sense.
The other $\Delta$-set are left empty.
Observe that, in particular, for all $q\ne q_0$ this yields
$$Minus^{[\gamma (i_0-1)+2]}(\widehat{q})
=Minus^{[\gamma (i_0-1)+1]}(\widehat{q})=q^{(i_0)}.$$
We then continue defining 
\begin{align*}
&\Delta ^{[\gamma (i_0-1)+2]}(Surplus(\widehat{q}_2))\\
&=
\powAst ( \{\Delta ^{[\gamma (i_0-1)+1]}(Surplus(\widehat{q}_1)) \} 
\cup\widehat{C_2} ^{[\gamma (i_0-1)+2]})\setminus
\{  \bigcup\widehat{C_2} ^{[\gamma (i_0-1)+2]}\},
\end{align*}
and all the argument used in the previous step can be repeated.

This procedure will prosecuted until the end of the cycle is reached, that is,
the node $C_{n+1}$.
At this step we introduce a slight modification in the construction of the $\Delta$-sets.
Namely, we have to restore the cardinality of $Minus(\widehat{q})$, 
which was pertubed moving 
$t_0$ from the Minus to the Surplus portion, in order to trigger off the pumping procedure.
Hence, pick an element $t_1$ inside 
$$\powAst (\{\Delta (Surplus^{[\gamma (i_0-1)+n+1]}(\widehat{q}_n)\} 
\cup\widehat{C_{n+1}} ^{[\gamma (i_0-1)+n+1]})\setminus
\{\bigcup\widehat{C_{n+1}} ^{[\gamma (i_0-1)+n+1]}\}.$$
Since we are assuming that
at each step the cycle can distribute at least 3 new elements, the set
\begin{align*}
&\Delta ^{[\gamma (i_0-1)+n+1]} (Surplus(\widehat{q}_0))\\
&=\powAst (\{\Delta (Surplus^{[\gamma (i_0-1)+n+1]}(\widehat{q}_n)\} 
\cup\widehat{C_{n+1}} ^{[\gamma (i_0-1)+n+1]})\setminus
\{\bigcup\widehat{C_{n+1}} ^{[\gamma (i_0-1)+n+1]}\}
\setminus \{t_1\}
\end{align*}
is certainly not empty.
Then define
$$ \Delta ^{[\gamma (i_0-1)+1]}(Minus(\widehat{q}_0))=\{t_1\}.$$
Notice that $t_1$ is unused, and so will be kept along the entire pumping
procedure of pumping, since
it lies in the Minus portion of $\widehat{q}_0$, 
which is untouched in this segment of the new formative process.
As before, the procedure can prosecute $\aleph _0$-times.

Since $q^{(\lambda)}=\bigcup_{\nu<\lambda}q^{(\nu)}$ for every $q\in
\Places$ and every limit ordinal $\lambda\leqslant \xi$, it is clear that
$\widehat{q}^{[\omega =\gamma (i_0)]}$ is equal to 
$\bigcup_{i\in{\mathbb N}}\widehat{q}^{[(i_0-1)+i]}$ for all $q\in
\Places$, consistently the Minus-Surplus partition is defined for the stage $\omega$.

By construction, for all $q\in\Places$ such that $q\neq q_0$
$(Minus^{\gamma(i_0)}(\widehat{q}))$ is equal to $q^{(i_0)}$ while
$(Minus^{[\gamma(i_0)]}(\widehat{q}_0))$ 
is equal to  $(q_0^{(i_0)}\setminus\{ t_0 \})\cup\{ t_1 \}$.

Our aim is to show that the transitive partitions 
$\Sigma_{i_0}$ and $\widehat{\Sigma}_{\gamma(i_0)}$ verify the
conditions to apply subsequently Lemma \ref{pasting} and Corollary \ref{upwards}, so proving the satisfiability of $\Phi$.

Concerning the application of Lemma~\ref{pasting},
we have to show properties (i), (vii), (viii), (x), and (a)-(c).
This is just a bookkeeping argument, and we
detail it in the Appendix.

Now the formative process $[\bullet]$ has copied the original one along 
the segment $[i_0,\ell]$. 
In order to apply
Lemma \ref{upwards}, we need to show that $\widehat{\Sigma}_{\gamma(\ell)}$ 
has the same target as $\Sigma_{\ell}$.
We simply observe that, if $q^{(\ell)}$ is a target of $\Gamma^{(\ell)}$,
there must
exist a step $i$ such that $\Gamma=A_{i}$ and $\Delta^{i}(q)\ne\emptyset$.
Since 
both the segment $[0,(i_0-1)]$ is equal to $[0,\gamma (i_0-1)]$,
and the segment $[i_0,\ell]$ is imitated by one application of Lemma \ref{pasting},
then 
$\Delta^{[\gamma (i)]}(\widehat{q})\ne\emptyset$ too.
On the other side, if $\Delta^{[\alpha]}(\widehat{q})\ne\emptyset$ for some $\alpha$, $q$
has to be a target of $\widehat{A}_{\alpha}$, so that we are done.

At this point Corollary \ref{upwards} applies, therefore 
all literals except those of $Finite$-type are satisfied.
Finally, the literals as $Finite(x)$ are satisfied as well.
Indeed, every block $q$ contained in $\Im_{{\mathcal M}}(x)$ 
lies in $\REDS$, and the formative process
$[\bullet]$ does not change size of such a block. Also, by hypothesis,
$\mbox{ for each variable } x$ such that $\neg Finite(x)\in\Phi,
\Im_{{\mathcal M}}(x)\cap\PLACES{\mathcal{C}}\mbox{ is not empty }$, 
and the blocks in the pumping cycle are infinitely increased during the pumping procedure.
Hence all of them are of
infinite size, as well as all the variables containing at least one of them.
This in turns implies that all $\neg Finite(x)\in\Phi$ are satisfied by the new model.
\qed

The above technique provides a valid tool to 
solve problems which require to build an infinite model.
In \cite{CU} it is shown that there is a computable function 
$f(n)$ such that, if a formula $\Phi$ in MLSSPF
is satisfiable, then there is an assignment rank bounded by 
$f(\vert{\mathcal X}_{\Phi}\vert)$ which satisfies a slight 
modification of the properties 
described in Theorem \ref{pumping1}. 
But then MLSSPF has the witness small property, and is
therefore decidable.
A similar argument it is used to prove the witness small property for MLSSPU.

\section{Open Problems}
\subsection{A Decidability Problem}
Even if all the problems related to the literals which force the infinity 
are treatable by the present approach, 
the decidability of MLSSP extended by the cartesian product binary 
operator [$x=y\times z$] is still an open question.
Observe that this language forces the infinity.
This problem is originally due to M.\ Davis, who proposed it as a set computable
version of the Tenth Hilbert Problem (see \cite{Mat}).

\subsection{ A Complexity Problem}
Decidability of MLSSP is NP-complete, therefore there is no hope to find
a polynomial time bound for our problems.
Nevertheless, the witness small model property furnishes 
double exponential decision algorithms. An exponential bound could be a good platform to 
perform polynomial time for special cases.

\section*{Appendix}
Here we exhibit a complete verification of the properties 
requested for the application of Lemma \ref{pasting}
within the proof of Theorem \ref{pumping1}.
\begin{itemize}
\item[(i)] First assume $q\ne q_{0}$. By construction, only  Surplus 
sides are increased along pumping procedure. 
Therefore $q^{(i_0)}=Minus^{[\gamma(i_0)]}\widehat{q}$. Otherwise, observe that
$Minus^{[\gamma(i_0)]}\widehat{q}=(q^{(i_0)}\setminus\{t_0\})\cup \{t_1\}$, hence
$\vert q^{(i_0)}\vert=\vert Minus^{[\gamma(i_0)]}\widehat{q}\vert$.

\item[(vii)] Observe that $\overline {\mathcal{C}}$ is 
composed of  green blocks only. Therefore, if $q\in\REDS$, by
hypothesis $q$ cannot belong to 
$\PLACES{\mathcal{C}}$, but the only blocks whose size is increased are inside 
$\PLACES{\mathcal{C}}$, hence 
$q^{(i_0)}=Minus^{[\gamma(i_0)]}\widehat{q}=\widehat{q}^{[\gamma(i_0)]}$.

\item[(viii)] Trivial.

\item[(x)] Assume $q_{0}\notin\Gamma$.
In this case, $Minus^{[\gamma(i_0)]}\widehat\Gamma=
Minus^{[\gamma(i_0-1)]}\widehat\Gamma$. Therefore, for all block $q$,
$$\powAst (Minus^{[\gamma (i_0-1)]}(\widehat{\Gamma}))\cap q^{[\gamma (i_0-1)+1]}
=\powAst (Minus^{[\gamma (i_0)]}(\widehat{\Gamma}))\cap q^{[\gamma (i_0-1)+1]}.$$
Along the pumping procedure, only the Surplus nodes are used. 
Since $\powAst$ of the Surplus nodes
are always disjoint from the Minus ones, we can prolongate the previous chain 
of equalities with
$$\powAst (Minus^{[\gamma (i_0)]}(\widehat{\Gamma}))\cap q^{[\gamma (i_0-1)+1]}=
\powAst (Minus^{[\gamma (i_0)]}(\widehat{\Gamma}))\cap q^{[\gamma (i_0)]}.$$
On the other hand, 
$$\powAst (\Gamma ^{(i_0-1)})\cap q^{(i_0)}
=\powAst (\Gamma ^{(i_0)})\cap q^{(i_0)}.$$
Finally, by construction,
$$\powAst (\Gamma ^{(i_0-1)})\cap q^{(i_0)}
=\powAst (Minus^{[\gamma (i_0-1)]}(\widehat{\Gamma}))\cap\widehat{q}^{[\gamma (i_0-1)+1]}.$$
In the other case, observe that $t_1$ is new at the step 
$\gamma(i_0-1)+n+1$. Thus everything created from $t_1$
cannot be inside any block $q$ before its distribution, neither in the segment 
$[\gamma (i_0-1)+n+1,\gamma (i_0)]$, for only the Surplus nodes are used, and $t_1$ is in the
Minus side of block $q_0$. This yields 
$$\powAst (Minus^{[\gamma (i_0-1)]}(\widehat{\Gamma}))\cap q^{[\gamma (i_0-1)+1]}
=\powAst (Minus^{[\gamma (i_0)]}(\widehat{\Gamma}))\cap q^{[\gamma (i_0-1)+1]}.$$
The prosecution of the argument follows exactly the one of the former case.

\item[(a)] If $q_{0}\in\Gamma$, the property trivially holds since $t_0$ is new at
the step $i_0$; therefore $\bigcup\Gamma$ cannot have been distributed at the stage $i_0$.
On the other hand $t_1$, which belongs to $Minus^{[\gamma (i_0-1)+n+1]}(\widehat{\Gamma}))$,
is new at the step $\gamma(i_0-1)+n+1$.
Hence $\bigcup Minus^{[\gamma (i_0-1)+n]}(\widehat{\Gamma})$ 
cannot have been distributed at the stage
$\gamma (i_0-1)+n]$. Again, the Minus nodes are unused along the pumping procedure, hence 
$\bigcup Minus^{[\gamma (i_0)]}(\widehat{\Gamma}))$ is not distributed 
at the limit step $\gamma (i_0)$ as well.
Conversely, if $q_{0}\notin\Gamma$, the result easily follow 
by standard arguments from the 
fact that the Minus portion of $\Gamma$ and the original $\Gamma$ are 
equal at the stage $\gamma (i_0)$, and 
the Minus nodes are unused along the pumping procedure.

\item[(b)] $Surplus^{[\gamma (i_0)]}(q)\ne\emptyset$ and $q\in\Gamma$,
therefore the node $\Gamma$
is changed along the pumping procedure. 
By construction, $\bigcup\Gamma$ is never distributed along pumping procedure, so 
$$\bigcup\widehat{\Gamma}^{[\gamma (i_0)]}
\in\powAst (\widehat{\Gamma}^{[\gamma (i_0)]})\setminus\bigcup_{q\in\Sigma}
\widehat{q}^{[\gamma (i_0)]}.$$

\item[(c)] Easily follows from the fact that  
after a grand event nothing changes in the formative process,
and from (ii) of Def.\ref{pumpchain}, 
which asserts that $GE(\PN{\PLACES{\mathcal{C}}}))\ge i_0$.
\end{itemize}


\end{document}